\documentclass[12pt]{article}
\usepackage{amssymb,bm,epsf,latexsym,graphicx}
\usepackage[english]{babel}
\def\bbbr{{\mathbb R}}
\def\bbbc{{\mathbb C}}

\def\bbbp{{\mathbb P}}
\def\bbbz{{\mathbb Z}}
\def\boldalpha{\bm{\alpha}}
\def\boldgamma{\bm{\gamma}}
\def\boldtheta{\bm{\theta}}
\def\boldtau{\bm{\tau}}

\def\boldrho{\bm{\rho}}
\def\boldmu{\bm{\mu}}
\def\boldsigma{\bm{\sigma}}

\def\mod#1{({\rm mod\ }#1)}

\def\v#1{{\bf#1}}

\def\D{\partial}

\def\Arg{{\rm Arg}}
\def\hypergeo#1#2#3#4#5{\,_{#1}F_{#2}\left(\left.{#3\atop #4}\right|#5\right)}
\def\is{\equiv}
\def\mod#1{({\rm mod}\ #1)}
\newtheorem{theorem}[subsection]{Theorem}

\newtheorem{lemma}[subsection]{Lemma}
\newtheorem{conjecture}[subsection]{Conjecture}
\newtheorem{corollary}[subsection]{Corollary}

\newtheorem{proposition}[subsection]{Proposition}

\newtheorem{assumption}[subsection]{Assumption}
\newenvironment{proof}{\medskip \noindent {\bf Proof}}
{\hfill\break\medskip\hfill$\Box$\medskip}

\title{Monodromy of A-hypergeometric functions}
\author{Frits Beukers}
\date{May 7, 2013}

\begin{document}
\maketitle

\abstract{Using Mellin-Barnes integrals we give a method to compute a
relevant subgroup of the monodromy group of an A-hypergeometric system
of differential equations. Presumably this group is the full monodromy
group of the system}

\parindent=0pt

\section{Introduction}
At the end of the 1980's Gel'fand, Kapranov and Zelevinsky,
in \cite{GKZ1}, \cite{GKZ2}, \cite{GKZ3}, defined a general
class of hypergeometric functions, encompassing the classical one-variable
hypergeometric functions, the Appell and Lauricella functions and Horn's functions.
They are called A-hypergeometric functions and they provide a beautiful and elegant
basis of a theory of hypergeometric functions in several variables. For
an introduction to the subject we refer the reader to \cite{stienstranotes}, \cite{beukersnotes}
or the book by Saito, Sturmfels and Takayama, \cite{SST}. We briefly recall the main facts.
Let $A\subset \bbbz^r$ be a finite set such that
\begin{enumerate}
\item The $\bbbz$-span of $A$ is $\bbbz^r$.
\item There exists a linear form $h$ such that $h(\v a)=1$ for
all $\v a\in A$.
\end{enumerate}
Let $\boldalpha=(\alpha_1,\ldots,\alpha_r)\in\bbbr^r$.
Denote $A=\{\v a_1,\ldots,\v a_N\}$ (with $N>r$). Writing the vectors $\v a_i$ in
column form we get the so-called A-matrix
$$A=\pmatrix{
a_{11} & a_{12} & \cdots & a_{1N}\cr
a_{21} & a_{22} & \cdots & a_{2N}\cr
\vdots &        &        & \vdots\cr
a_{r1} & a_{r2} & \cdots & a_{rN}\cr
}
$$
For $i=1,2,\ldots,r$ consider the first order differential operators
$$Z_i=a_{i1}v_1\D_1+a_{i2}v_2\D_2+\cdots+a_{iN}v_N\D_N$$
where $\D_j={\D\over \D v_j}$ for all $j$.

Let
$$L=\{(l_1,\ldots,l_N)\in\bbbz^N|\ l_1\v a_1+l_2\v a_2+\cdots+l_N\v a_N=\v 0\}$$
be the lattice of integer relations between the elements of $A$. For every $\v l\in L$
we define the so-called box-operator
$$\Box_{\v l}=\prod_{l_i>0}\D_i^{l_i}-\prod_{l_i<0}\D_i^{-l_i}$$
The system of differential equations
\begin{eqnarray*}
(Z_i-\alpha_i)\Phi&=&0\quad (i=1,\ldots,r)\\
\Box_{\v l}\Phi&=&0\quad \v l\in L
\end{eqnarray*}
is known as the system of {\it A-hypergeometric differential equations} and we denote it
by $H_A(\boldalpha)$.
It turns out that in general the solution space of $H_A(\boldalpha)$ is finite dimensional
with dimension equal to the volume of the convex hull $Q(A)$ of $A$. In order to be more precise we
have to introduce $C(A)$, the cone generated by the $\bbbr_{\ge0}$-linear combinations of
$\v a_1,\ldots,\v a_N$.
We say that an A-hypergeometric system is {\it non-resonant} if the boundary of $C(A)$
has empty intersection with the
shifted lattice $\boldalpha+\bbbz^r$. We have the following theorem.

\begin{theorem}[GKZ, Adolphson]\label{dimension}
Suppose either one of the following conditions holds,
\begin{enumerate}
\item the toric ideal $I_A$ in $\bbbc[\D_1,\ldots,\D_N]$ generated by the box
operators has the Cohen-Macaulay property.
\item The system $H_A(\boldalpha)$ is non-resonant.
\end{enumerate}
Then the rank of $H_A(\boldalpha)$ is finite and equals
the volume of the convex hull $Q(A)$ of the points of $A$. The volume is normalized so that a minimal
$(r-1)$-simplex with integer vertices in $h(\v x)=1$ has volume $1$.
\end{theorem}

Theorem \ref{dimension} is proven in \cite{GKZ3},
(corrected in \cite{GKZ5}) and \cite[Corollary 5.20]{adolphson}.

Among the many papers written on A-hypergeometric equations there are very few papers
dealing with the monodromy group of these systems in general. In the case of one-variable hypergeometric
functions there is the paper by Beukers and Heckman, \cite{beukersheckman}, which give an characterisation of monodromy groups
as complex reflection groups. There is also a classical method to compute monodromy with respect
to an explicit basis of functions using so-called Mellin-Barnes integrals, see F.C.Smith, \cite{smith}. In the special case of the fourth order equation with symplectic monodromy
there is a detailed calculation in \cite{ChenYangYui}.
A recent paper by Golyshev and Mellit, \cite{golyshevmellit}, deals with the same problem using Fourier-transforms of $\Gamma$-products. A recent paper by K.Mimachi \cite{mimachi}
uses computation of twisted cycle intersection.

For the two-variable Appell system $F_1$ and Lauricella's $F_D$, monodromy follows from the work by
E.Picard \cite{picard}, T.Terada \cite{terada} and Deligne-Mostow \cite{delignemostow}.
In T.Sasaki's paper \cite{sasaki} we find explicit monodromy generators for Appell $F_1$.
They all use the fact that Lauricella functions of $F_D$-type can be written as one-dimensional twisted period integrals and monodromy is a representation of the pure braid group on
$n+1$ strands (where $n$ is the number of variables).

The two-variable Appell $F_2$ has been considered explicitly by M.Kato, \cite{kato}. The
Appell $F_3$ system has the same $A$-set as Appell $F_2$, and therefore gives nothing new.
Finally, the Appell system $F_4$ has been
considered completely explicitely by K.Takano \cite{takano} and later Kaneko, \cite{kaneko}. In
Haraoka, Ueno \cite{haraokaueno} we find some rigidity considerations on the monodromy of $F_4$.
In the paper \cite{matsumotoyoshida} by K.Matsumoto and M.Yoshida, the authors provide generators
for the monodromy of Lauricella $F_A$.

Finally, the complete monodromy of the Aomoto-Gel'fand system $E(3,6)$ has been determined
by K.Matsumoto, T.Sasaki, N.Takayama, M.Yoshida in \cite{masatayo1} and further properties in
\cite{masatayo2}. See also M.Yoshida's book
'Hypergeometric Functions, my Love', \cite{yoshida}.

In essentially all of the above studies the monodromy is computed by studying the behaviour of
Euler integrals for hypergeometric functions under analytic continuation and corresponding
deformation of the contours of integration. For this, knowledge of the fundamental group
of the complement of the singular locus of the system of equations is required.
It is the purpose of the present paper to avoid these geometric difficulties as long as
possible and compute monodromy groups of A-hypergeometric systems
by methods which are combinatorial in nature. We do this by starting with local monodromy groups
which arise from series expansions of solutions of $H_A(\boldalpha)$. It is well known that such
local expansions correspond one-to-one with regular triangulations of $A$. This is a discovery
by Gel'fand, Kapranov and Zelevinsky that we shall explain in Section \ref{powerseries}.
The local monodromy groups have to be glued together to build a global monodromy group. This
glue is provided by multidimensional Mellin-Barnes integrals as defined in Section \ref{mellinbarnes}.
Unfortunately, the Mellin-Barnes integrals do not always provide a basis of solutions.
But if they do (Assumption \ref{mellinbasis}),
the construction of the global group generated by the local contributions is completely combinatorial. In Section \ref{algorithm} we give a practical recipe for the calculation
of these matrices. This algorithm is based on the theoretical considerations in the preceding
sections.

Although in a good number of classical cases the method described in Section \ref{algorithm}
works very well, it is not always garantueed to work. There are two potential obstacles:
\begin{enumerate}
\item There does not always exist a basis of Mellin-Barnes solutions. In many classical
cases such a basis exists. For example, two variable Appell, Horn and higher Lauricella $F_A,F_B,F_D$. But, on the other hand, in case of Lauricella $F_C$ and many other Aomoto-systems
such a basis does not seem to exist. This is a subject of further investigation.
\item The group we calculate is the subgroup generated by the contributions of local
monodromies at different points, modulo scalars. Let us call this group $lMon$.
It is not clear if
this group equals the complete monodromy group modulo scalars, which we denote $Mon$.

The reason we consider the monodromy group modulo scalars is that
for the A-hypergeometric system and their classical counterparts these groups are the same.
An explanation for this can be found at the end of Section \ref{powerseries}.
\end{enumerate}

The groups we calculate are determined with respect to
a basis of solutions in Mellin-Barnes integral form. In the case of one variable $_{n+1}F_n$
they turn out to coincide with the matrices found in \cite{beukersheckman} (see
Section \ref{example3F2}).
If one would like to calculate monodromy matrices with respect to explicit bases of
local power series expansions
one would have to find an explicit calculation of a Mellin-Barnes
integral as a linear combination of power series solutions. This is a tedious task
which we like to carry out in a forthcoming paper. We remark that such a calculation has been
carried out in the so-called confluent case (i.e. $A$ does not lie in translated hyperplane) by
O.N.Zhdanov and A.K.Tsikh, \cite{zhdanovtsikh}. In her PhD-thesis from 2009, Lisa Nilsson \cite{nilsson}
introduced (non-confluent) A-hypergeometric functions in terms of Mellin-Barnes integrals and initiated their study. It is these integrals that we shall use.

In the remainder of this paper we assume that {\it the shifted lattice $\boldalpha+\bbbz^r$ has empty
intersection with any hyperplane spanned by $r-1$ independent elements of $A$}. In that case we say that
the system is {\it totally non-resonant}. Note that this is stronger than just non-resonance where
only the faces of the cone spanned by the elements of $A$ are involved. Non-resonance
(and a fortiori total non-resonance) ensures that our system is irreducible,
see for example \cite[Thm 2.11]{GKZ4} or
\cite{beukersirreducible} for a slightly more elementary proof. Non-resonance also implies that A-hypergeometric
systems whose parameter vectors are the same modulo $\bbbz^r$ have isomorphic monodromy,
see \cite[Thm 2.1]{beukersirreducible}, which is actually a theorem due to B.Dwork. Total non-resonance implies
T-nonresonance for every triangulation $T$ in the terminology of Gel'fand, Kapranov and Zelevinsky.
In particular this implies that the local solution expansions will not contain logarithms. So
local monodromy representations act by characters. We prefer to leave the case of logarithmic local
solutions for a later occasion.

The starting data of our computation will not be the set $A$ and parameter vector $\boldalpha$,
but rather a dual version as follows. Let $d=N-r$. This will be the number of variables in
the classical counterpart of the A-hypergeometric system (number of essential variables,
e.g. $d=2$ in the Appell cases).
Choose a $\bbbz$-basis for the lattice $L$, which has rank $d$, and write
the basis elements as rows of a $d\times N$ matrix $B$. In the literature the transpose
of $B$ is often called a
{\it Gale dual} of $A$, we simply call $B$ a {\it B-matrix}. The matrix $B$ has the property
that it has maximal rank $d$, the $\bbbz$-span of the columns is $\bbbz^d$ and $A.B^t$ is the zero matrix. We denote the columns of $B$ by $\v b_j,\ j=1,\ldots,N$.
Then the space $L\otimes\bbbr\subset\bbbr^N$ is parametrized by the $N$-tuple
$(\v b_1\cdot\v s,\ldots,\v b_N\cdot\v s)$
with $\v s=(s_1,\ldots,s_d)\in\bbbr^d$ as parameters.
In our computations we take a B-matrix as starting data and instead of the parameters $\alpha$
we choose $\boldgamma=(\gamma_1,\ldots,\gamma_N)\in\bbbr^N$ such that
$\gamma_1\v a_1+\cdots+\gamma_N\v a_N=\boldalpha$. Notice that there is some ambiguity in
the choice of $\boldgamma$ which we will fix later. The reason we take $B$ and $\boldgamma$
as starting data is that they are easily read off from the classical power series expansions.
In the next section we see an example of this.

\section{Power series solutions}\label{powerseries}
Consider the system $H_A(\boldalpha)$ and a formal solution
$$\Phi_{\boldgamma}=\sum_{\v l\in L}{v_1^{l_1+\gamma_1}\cdots v_N^{l_N+\gamma_N}
\over\Gamma(l_1+\gamma_1+1)\cdots\Gamma(l_N+\gamma_N+1)}$$
where $\boldgamma$ is chosen such that $\boldalpha=\gamma_1\v a_1+\cdots+\gamma_N\v a_N$.
This expansion was introduced in \cite{GKZ3}. It is a Laurent series multiplied by generally non-integral
powers of the variables $v_i$. We call such a series a twisted Laurent series.
As is well-known we have a freedom of choice in $\boldgamma$ by shifts over
$L\otimes\bbbr$.
We shall use this freedom in the following way where, again, we denote the columns
of the B-matrix by $\v b_i$.
Choose a subset $I\subset\{1,2,\ldots,N\}$ with $|I|=d=N-r$ such that $\v b_i$
with $i\in I$ are linearly independent.
It is known that $|\det(\v b_i)_{i\in I}|=|\det(\v a_j)_{j\not\in I}|$ and we denote
this quantity by $\Delta_{I}$. Choose $\boldgamma$ such that
$\gamma_i\in\bbbz$ for all $i\in I$. There are precisely
$\Delta_I$ such choices for $\boldgamma$ which are distinct modulo $L\otimes\bbbr$.
The series $\Phi$ now reads
$$\Phi_{\boldgamma}=\sum_{\v l\in L}\prod_{i\in I}{v_i^{l_i+\gamma_i}\over\Gamma(l_i+\gamma_i+1)}
\times\prod_{j\not\in I}{v_j^{l_j+\gamma_i}\over\Gamma(l_j+\gamma_j+1)}$$
Since $\gamma_i\in\bbbz$ for all $i\in I$ this summation extends over all $\v l$
with $\gamma_i+l_i\ge0$ for all $i\in I$. The other terms vanish because
$1/\Gamma(\gamma_i+l_i+1)=0$ whenever $\gamma_i+l_i$ is a negative integer for some $i\in I$.
Hence $\Phi_{\boldgamma}$ is now a twisted power series. Let us fix such a choice of $\boldgamma$.
It is not hard to see that a set of series $\Phi_{\boldgamma}$ with $\boldgamma$-values
which are distinct modulo $L\otimes\bbbr$, is linearly independent over $\bbbc$.

Let now $\rho_1,\ldots,\rho_N$ be any $N$-tuple with the property that $\rho_1l_1+\cdots+\rho_Nl_N>0$
for any non-zero $\v l\in L$ with $l_i\ge 0$ for all $i\in I$.
Then, according to Theorem \cite[Proposition 16.2]{beukersnotes}
or \cite[Section 3.3,3.4]{stienstranotes}, the series $\Phi_{\boldgamma}$ converges for all $v_1,\ldots,v_N$ with
$\forall i:|v_i|=t^{\rho_i}$ and $t\in\bbbr_{>0}$ sufficiently small.
We call such an $N$-tuple $\rho_1,\ldots,\rho_N$ a {\it convergence direction} of $\Phi_{\boldgamma}$.

There is one important assumption we need in order to make this approach work. Namely
the garantee that none of the arguments $\gamma_j+l_j$ is a negative integer when $j\not\in I$.
Otherwise we might even end up with a trivial series solution. Notice that
$$\boldalpha=\sum_{j=1}^N\gamma_j\v a_j\is\sum_{j\not\in I}\gamma_j\v a_j\mod{\bbbz^r}$$
So if $\gamma_j\in\bbbz$ for some $j\not\in I$, the point $\boldalpha$ lies modulo $\bbbz^r$ in
a space spanned by the $r-1$ remaining vectors $\v a_i$. Under the assumption of total non-resonance
this situation cannot occur.

From now on we assume that $H_A(\boldalpha)$ is totally non-resonant. We denote the set of all sets
$I$ such that $\Delta_I=|\det(\v b_i)_{i\in I}|\ne0$ by ${\cal I}$. When $s=\Delta_I>1$ we
must take $s$ copies of $I$ in this list.
To each $I\in{\cal I}$ there corresponds a choice of $\boldgamma$ and we see to
it that all these choices are distinct modulo $L\otimes\bbbr$. So to an index set $I$ which occurs $s$ times there
correspond $s$ choices of $\boldgamma$ that are distinct modulo $L\otimes\bbbr$. The corresponding powerseries solutions are
denoted by $\Phi_I$.

Choose $I\in{\cal I}$ and a convergence direction
$(\rho_1,\ldots,\rho_N)$ such that $\rho_1l_1+\cdots+\rho_Nl_N>0$
for all non-zero $\v l\in L$ with $\forall i\in I:l_i\ge0$. Note that if
$(\rho_1,\ldots,\rho_N)$ is a convergence direction, then after adding an
element of the $\bbbr$-row span of $A$, it is still a convergence condition
since $A\v x=\v 0$ for all $\v x\in L\otimes\bbbr$. Consider the element
$\boldrho=\sum_{i=1}^N\rho_i\v b_i$ in the column span of $B$. By shifting over
the row span of $A$ we can see to it that $\rho_i=0$ for all $i\not\in I$.
Hence $\boldrho=\sum_{i\in I}\rho_i\v b_i$. The convergence condition is now
equivalent to saying that the new $\rho_i$ are positive. So if we denote
$$\v b_I=\left\{\left.\sum_{i\in I}\lambda_i\v b_i\right|\lambda_i>0\right\},$$
the convergence condition can be restated as $\boldrho\in\v b_I$.
By a slight abuse of language we call the vector $\boldrho$ also a convergence
direction.

Conversely, fix an element $\boldrho$ in the span of all $\v b_i$ which does not lie on the boundary of
any $\v b_I$. Define ${\cal I}_{\boldrho}=\{I|\rho\in\v b_I\}$.
Then, by the theory of Gel'fand, Kapranov and Zelevinsky the powerseries $\Phi_I$
with $I\in{\cal I}_{\boldrho}$ form a basis of solutions with a common open region of convergence.
We call such a set a basis of local solutions of $H_A(\boldalpha)$.
It also follows from the theory that the sets ${\cal I}_{\boldrho}$ are in one-to-one correspondence with
the regular triangulations of the set $A$. This correspondence is given by
${\cal I}_{\boldrho}\mapsto\{I^c|I\in {\cal I}_{\boldrho}\}$.

The intersections of the simplicial cones $\v b_I$ define a subdivision of
$\bbbr^d$ into open convex polyhedral cones whose
closure of the union is $\bbbr^d$. This is a polyhedral fan which is called the {\it secondary fan}.
The open cones in the secondary fan are in one-to-one correspondence with the bases of local series
solutions. As an example let us take the system Appell $F_2$. The standard Appell $F_2$-series
reads
$$F_2(\alpha,\beta,\beta',\gamma,\gamma',x,y)=\sum_{m,n\ge0}
{(\alpha)_{m+n}(\beta)_m(\beta')_n\over m!n!(\gamma)_m(\gamma')_n}\ x^my^n.\label{F2expansion1}
$$
We hope no confusion arise with the existing notations $\boldalpha,\boldgamma$.
Using the identity
$$\Gamma(z)\Gamma(1-z)=\pi/\sin\pi z$$
we see that the series is proportional to

$$\sum_{m,n\ge0}{\scriptstyle{x^my^n}\over\scriptstyle{ \Gamma(-\alpha-m-n+1)\Gamma(-\beta-m+1)\Gamma(-\beta'-n+1)
\Gamma(\gamma+m)\Gamma(\gamma'+n)\Gamma(m+1)\Gamma(n+1)}}.\label{F2expansion2}$$

The basisvectors $(-1, -1, 0, 1, 0, 1, 0)$ and $(-1, 0, -1, 0, 1, 0, 1)$ of $L$
are given to us naturally because these are the coefficient vectors of $m$ and
$n$ respectively in the $\Gamma$-factors of the expansion just given. This follows
from the shape of the canonical solution $\Phi_{\boldgamma}$. So our B-matrix reads
$$B=\pmatrix{-1 & -1 & 0 & 1 & 0 & 1 & 0\cr -1 & 0 & -1 & 0 & 1 & 0 & 1\cr}.$$
A parameter vector $\boldgamma$ can also be read off from the $\Gamma$-expansion, namely
$$\boldgamma=(-\alpha,-\beta,-\beta',\gamma-1,\gamma'-1,0,0).$$
The column vectors of $B$ are depicted here,

\centerline{\includegraphics[height=5cm]{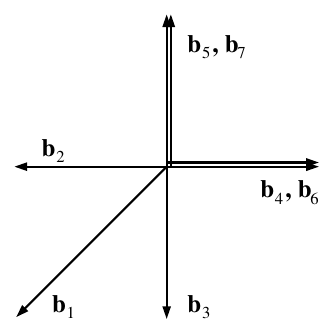}}

For example, consider the vector $(-0.5,1)$ in this picture. We see that it is
contained in the positive cones of the following pairs: $\{\v b_2,\v b_5\},
\{\v b_2,\v b_7\},\{\v b_1,\v b_5\},\{\v b_1,\v b_7\}$.
Taking the complimentary sets of indices of each pair we get
$$\{1,3,4,6,7\},\{1,3,4,5,6\},\{2,3,4,5,6\},\{2,3,4,5,6\}.$$
These form the index sets of the simplices of a triangulation of the set $A$.
Take the alternative parameter vector
$$\boldgamma=(\gamma'-\beta-\alpha-1,0,\gamma'-1-\beta',\gamma-\beta-1,0,-\beta,1-\gamma'),$$
which differs from the original choice by an element of $L\otimes\bbbr$.
The coordinates on positions 2,5 are made zero, corresponding to the choice
$\v b_2,\v b_5$. The formal solution for this new parameter vector reads

$$\sum_{m,n}{\scriptstyle{v_1^{\gamma'-\beta-\alpha-1-m-n}v_2^{-m}v_3^{\gamma'-\beta'-1-n}
v_4^{\gamma-\beta-1+m}v_5^{n}v_6^{-\beta+m}v_7^{1-\gamma'+n}}
\over\scriptstyle{\Gamma(\gamma'-\beta-\alpha-m-n)\Gamma(-m+1)\Gamma(\gamma'-\beta'-n)
\Gamma(\gamma-\beta+m)\Gamma(n+1)\Gamma(-\beta+m)\Gamma(2-\gamma'+n)}}.$$

Since $1/\Gamma(n+1)=0$ when $n<0$ and $1/\Gamma(-m+1)=0$ when $m>0$ we see that our summation
runs over $m\le0$ and $n\ge0$. Replace $m$ by $-m$ to get
$\scriptstyle{v_1^{\gamma'-\beta-\alpha-1}v_3^{\gamma'-\beta'-1}
v_4^{\gamma-\beta-1}v_6^{-\beta}v_7^{1-\gamma'}}$ times

$$\sum_{m,n\ge0}{\scriptstyle{
(v_1v_2v_4^{-1}v_6^{-1})^m(v_5v_7v_1^{-1}v_3^{-1})^n}
\over\scriptstyle{\Gamma(\gamma'-\beta-\alpha+m-n)\Gamma(m+1)\Gamma(\gamma'-\beta'-n)
\Gamma(\gamma-\beta-m)\Gamma(n+1)\Gamma(-\beta-m)
\Gamma(2-\gamma'+n)}},$$

a twisted powerseries in $v_1v_2v_4^{-1}v_6^{-1}$ and $v_5v_7v_1^{-1}v_6^{-1}$,
which has $(-0.5,1)$ as convergence direction. In the same way we can construct three other
(twisted) powerseries expansions and thus obtain a basis of local powerseries solutions of
our system.

To return to our general story, suppose
we have a basis of local series solutions
and suppose they converge in an open neighbourhood of a set
defined by $|v_1|=r_1,\ldots,|v_N|=r_N$ with $r_i>0$ for all $i$.
Let $\boldgamma^{(1)},\ldots,\boldgamma^{(D)}$ be the
choices of $\boldgamma$ used in the construction of the local series solutions $\Phi_1,\ldots,\Phi_D$.
Given any $N$-tuple
of integers $n_1,\ldots,n_N$ we define the closed path $c(n_1,\ldots,n_N)$ by
$$(r_1e^{2\pi i n_1t},\ldots,r_Ne^{2\pi i n_Nt}),\quad t\in[0,1].$$
Taking $(r_1,\ldots,r_N)$ as a base point, the series $\Phi_j$ changes into
$$\exp(2\pi i (n_1\gamma_1^{(j)}+\cdots+n_N\gamma_N^{(j)}))=\exp(2\pi i
\ \v n\cdot\boldgamma^{(j)})$$
times $\Phi_j$ after analytic continuation along $c(n_1,\ldots,n_N)$.
The group of these substitutions is called a {\it local monodromy
group} with respect to $\Phi_1,\ldots,\Phi_D$. It is important to note that it
is generated by $d$ elements and a group of scalar elements generated by
the scalars $\exp(2\pi i\alpha_j),\ j=1,\ldots,r$. The explanation is
as follows. The $r$ rows of the matrix $A$ span a lattice. The basis of rows
can be completed to a basis of $\bbbz^N$ by $d=N-r$ extra integer vectors
$\v n_1,\ldots,\v n_d$, say. So all $\v n\in\bbbz^N$ are $\bbbz$-linear
combinations of these vectors. Suppose that $\v n$ is the $j$-th row of $A$.
Since $\v n\cdot\boldgamma^{(i)}=\alpha_j$ for all $i=1,\ldots,D$, the local monodromy
transformation is the scalar element given by $\exp(2\pi i\alpha_j)$.
The remaining vectors $\v n_1,\ldots,\v n_d$ provide us with $d$ generators.
In our implementation in Section \ref{algorithm} we shall make a sensible choice
for these generators.

The fact that we have $d$ generators corresponds to the fact that the number of
essential variables (modulo homogeneities) is $d$. The scalar group is simply the
difference between the full A-hypergeometric system and its classical counterpart.
In this paper we shall adopt the convention that we compute the monodromy modulo
scalars. Hence it suffices to compute the $d$ generators of the local monodromies.

\section{Mellin-Barnes integrals}\label{mellinbarnes}
Let notations be as in the previous sections. Consider the parametrization
of $L\otimes\bbbr$ by the $N$-tuple
$(\v b_1\cdot\v s,\ldots,\v b_N\cdot\v s)$
with $\v s=(s_1,\ldots,s_d)\in\bbbr^d$ as parameters.
Choose $\boldsigma=(\sigma_1,\ldots,\sigma_d)\in\bbbr^d$.

We now complexify the parameters $s_i$ and consider
the integral
$$M(v_1,\ldots,v_N)=
\int_{\boldsigma+\sqrt{-1}\ \bbbr^d}\prod_{i=1}^N \Gamma(-\gamma_i-\v b_i\cdot\v s)
v_i^{\gamma_i+\v b_i\cdot\v s}d\v s
\eqno{(MB)}$$
where $d\v s=ds_1\wedge\cdots\wedge ds_d$ and the integration takes place over
$-\infty<{\rm Im}(s_i)<\infty$  and ${\rm Re}(s_i)=\sigma_i$ for $i=1,\ldots,d$.
This is an example of a so-called {\it Mellin-Barnes
integral}. It will be crucial in the determination of the monodromy of A-hypergeometric
systems.

We prove the following Theorem,
\begin{theorem}\label{mbsolution}
Assume that $\gamma_i<-\v b_i\cdot\boldsigma$ for $i=1,2,\ldots,N$. Then the Mellin-Barnes integral
$M(v_1,\ldots,v_N)$ satisfies the set of A-hypergeometric equations $H_A(\boldalpha)$.
\end{theorem}

This will be done under the assumption that the Mellin-Barnes integral converges absolutely.
We come to the matter of convergence in the next section.

\begin{proof}: The Mellin-Barnes integral clearly has the property
$$M(\v t^{\v a_1}v_1,\ldots,\v t^{\v a_N}v_N)
=\v t^{\boldalpha}M(v_1,\ldots,v_N)$$ for all $\v t\in(\bbbc^*)^N$.
So $M(\v v)$ satisfies
the hypergeometric homogeneity equations.

Now let $\lambda\in L$ and put $\lambda=\lambda_+-\lambda_-$ where $\lambda_{\pm}$ have
non-negative coefficients and disjoint support. Define
$|\lambda|=\sum_{i=1}^N|\lambda_i|$. Then
\begin{eqnarray*}
&&\Box_{\lambda}M(v_1,\ldots,v_N)\\
&=&(-1)^{|\lambda|/2}\int_{\boldsigma+\sqrt{-1} \bbbr^d}
\prod_{i=1}^N\Gamma(-\gamma_i-\v b_i\cdot\v s+\lambda_{+,i})
v_i^{\gamma_i+\v b_i\cdot\v s-\lambda_{+,i}}d\v s\\
&-&(-1)^{|\lambda|/2}\int_{\boldsigma+\sqrt{-1}\bbbr^d}
\prod_{i=1}^N\Gamma(-\gamma_i-\v b_i\cdot\v s+\lambda_{-,i})
v_i^{\gamma_i+\v b_i\cdot\v s-\lambda_{-,i}}d\v s.
\end{eqnarray*}
Choose $\v s_{\lambda}$ such that $\v b_i\cdot\v s_{\lambda}=\lambda_i$ for $i=1,\ldots,N$. Then the second
integral is actually integration over $\v s_{\lambda}+\sqrt{-1}\bbbr^d$ of the integrand of the first integral.

Because of the assumption $\gamma_i<-\v b_i\cdot\boldsigma$ we see that
$-\gamma_i-t\v b_i\cdot(\v s_{\lambda}+\boldgamma)+\lambda_{+,i}>0$ for
all $t\in[0,1]$ and $i=1,\ldots,N$. Hence the $d+1$-dimensional domain
$\{t(\v s_{\lambda}+\boldgamma)+i\bbbr^d|0\le t\le1\}$ does not contain
any poles of the integrand and a homotopy argument gives that the two integrals are equal and cancel.
\end{proof}

Not all systems $H_A(\boldalpha)$ allow a choice of negative $\gamma_i$. However, under the assumption of non-resonance
it is known that two hypergeometric systems $H_A(\boldalpha)$ and $H_A(\boldalpha')$ have the same monodromy if
$\boldalpha-\boldalpha'\in\bbbz^r$ (see \cite[Theorem 2.1]{beukersirreducible}).
We call such systems contiguous.
Thus we can always replace an irreducible A-hypergeometric system by a contiguous one
which does allow
a choice of $\gamma_i<-\v b_i\cdot\boldsigma$ for all $i$. In concrete cases we
can also play with the value of $\boldsigma$.
From now on we make this assumption, i.e all our Mellin-Barnes integrals
are solution of an A-hypergeometric system.
Of course there is also the question whether or not $M(v_1,\ldots,v_N)$ is a trivial function. By
Proposition \ref{mellinbasis} we will find that it is non-trivial.

\section{Convergence of the Mellin-Barnes integral}
We find from \cite{AAR} the following estimate. Suppose $s=a+bi$ with $a_1<a<a_2$ and $|b|\to\infty$.
Then
$$|\Gamma(a+bi)|= \sqrt{2\pi}|b|^{a-1/2}e^{-\pi|b|/2}[1+O(1/|b|)].$$
Notice also that for any $v\in\bbbc^*$ we have $|v^{a+bi}|=|v|^ae^{-b\arg(v)}$.
Write $s_j=\sigma_j+i\tau_j$ for $j=1,\ldots,N-r$.
Let us denote $\theta_j=\arg(v_j)$ and $l_j(\boldtau)=l_j(\tau_1,\ldots,\tau_d)$.
The integrand in the Mellin-Barnes integral can now be estimated by
$$\left|\prod_{i=1}^N\Gamma(-\gamma_i-\v b_i\cdot\v s)) v_i^{\gamma_i+\v b_i\cdot\v s}\right|
\le c_1\max_j|\tau_j|^{c_2}\exp\left(-\sum_{j=1}^N\pi|\v b_j\cdot\boldtau|/2-\theta_j\v b_j\cdot\boldtau\right)$$
where $c_1,c_2$ are positive numbers depending only on $\gamma_j,v_j,\sigma_j$.
In order to ensure convergence of the integral we must have that
$$\sum_{j=1}^N\pi|\v b_j\cdot\boldtau|/2+\theta_j\v b_j\cdot\boldtau >0
\eqno{{\rm (C)}}\label{convergence}$$
for every non-zero $\boldtau\in\bbbr^d$.
We apply the following Lemma.

\begin{lemma}\label{convcondition}
Given $N$ vectors $\v p_1,\ldots,\v p_N$ in $\bbbr^d$. Suppose they have
rank $d$. Let $\v q\in\bbbr^d$. Then the following statements are equivalent,
\begin{itemize}
\item[i)] For all non-zero $\v x\in\bbbr^d$:
$$|\v q\cdot\v x|<\sum_{i=1}^N|\v p_i\cdot\v x|$$
\item[ii)] There exist $\mu_1,\ldots,\mu_N$ with $-1<\mu_i<1$ such that
$$\v q=\mu_1\v p_1+\cdots+\mu_N\v p_N.$$
\end{itemize}
\end{lemma}

\begin{proof}: First suppose that $\v q=\mu_1\v p_1+\cdots+\mu_N\v p_N$. Then, for all
non-zero $\v x\in\bbbr^d$,
\begin{eqnarray*}
|\v q\cdot\v x|&=&|\sum_{i=1}^N\mu_i \v p_i\cdot\v x|\\
&\le& \sum_{i=1}^N|\mu_i||\v p_i\cdot\v x|<\sum_{i=1}^N|\v p_i\cdot\v x|
\end{eqnarray*}

To show the converse statement consider the set
$$V=\left\{\left.\sum_{i=1}^N \mu_i\v p_i\ \right|\ -1<\mu_i<1\right\}.$$
This is a convex set. Suppose $\v q\not\in V$. Then there exists a linear form $h$ such that
$h(\v q)>h(\v p)$ for all $\v p\in V$. In other words, there exists $\v x\in\bbbr^d$ such that
$\v q\cdot\v x>\sum_{i=1}^N\lambda_i\v p_i\cdot\v x$ for all
$-1<\lambda_i<1$. In particular,
$$|\v q\cdot\v x|\ge\sum_{i=1}^N|\v p_i\cdot\v x|$$
contradicting our assumption. Hence $\v q\in V$.
\end{proof}

Application of Lemma \ref{convcondition} with $\v q=\sum_{j=1}^N\theta_j\v b_j$
and $\v p_j=\pi\v b_j/2$ to inequality (C) on page \pageref{convergence}
yields the following criterion.

\begin{corollary}\label{mbconverge}
Let notations be as above. Then the Mellin-Barnes integral converges absolutely
if there exist $\mu_i\in(-1,1)$ such that
$$\sum_{i=1}^N{\theta_i\over2\pi}\v b_i={1\over4}\sum_{i=1}^N\mu_i\v b_i.$$
\end{corollary}

Let us define
$Z_B=\{{1\over4}\sum_{i=1}^N\mu_i\v b_i|\mu_i\in(-1,1)\}$. This
is a so-called zonotope in $d$-dimensional space. The convergence condition for
the Mellin-Barnes integral now reads
$$\sum_{i=1}^N{\theta_i\over2\pi}\v b_i\in Z_B.$$

As an example let us again take the system Appell $F_2$. Recall
that our B-matrix reads
$$B=\pmatrix{-1 & -1 & 0 & 1 & 0 & 1 & 0\cr -1 & 0 & -1 & 0 & 1 & 0 & 1\cr}^t.$$
As before, a parameter vector $\boldgamma$ can also be read off from the
$\Gamma$-expansion, namely
$\boldgamma=(-\alpha,-\beta,-\beta',\gamma-1,\gamma'-1,0,0)$.
When $\alpha,\beta,\beta'>0$ and $\gamma,\gamma'<1$ we see that $\boldgamma$ has negative
components, except for the last two. By making a suitable choice for $\boldsigma\in\bbbr^2$
we can see to it that $\gamma_i<-\v b_i\cdot\boldgamma$ for all $i$.
Thus the corresponding Mellin-Barnes integral is indeed a solution of the $F_2$-system.
The zonotope $Z_B$ can be pictured as

\centerline{\includegraphics[height=5cm]{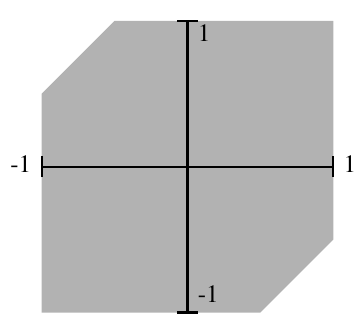}}

and the convergence condition reads
$${1\over2\pi}(-\theta_1-\theta_2+\theta_4+\theta_6,
-\theta_1-\theta_3+\theta_5+\theta_7)\in Z_B.$$

Note that the four points $(\pm 1/2,\pm 1/2)$ are contained
in $Z_B$. They correspond to the arguments
$\boldtheta=2\pi (0,0,0,0,0,\pm1/2,\pm1/2)$.
These argument choices
represent the same point in $v_1,\dots,v_7$\,-space. Hence we have four Mellin-Barnes
solutions of the system $F_2$ around one point. According to Proposition \ref{mellinindependent} these
integrals are linearly independent, and hence form a basis of local solutions of the $F_2$-system.
We say that we have a {\it Mellin-Barnes basis} of solutions.

\begin{proposition}\label{mellinindependent}
Let $\v v_0=(v^{(0)}_1,\ldots,v^{(0)}_N)\in(\bbbc^*)^N$ and let $\Theta$ be a finite set of
$N$-tuples $\boldtheta=(\theta_1,\ldots,\theta_N)$
such that $v^{(0)}_j=|v^{(0)}_j|\exp(i\theta_j)$ and the sums ${1\over2\pi}\sum_{j=1}^d\theta_j\v b_j$ are distinct
elements of $Z_B$. To each $\boldtheta\in \Theta$ denote
the corresponding determination of the Mellin-Barnes integral in the neighbourhood of $\v v_0$ by $M_{\boldtheta}$.

Then the functions $M_{\boldtheta}$ are linearly independent over $\bbbc$.
\end{proposition}

The proof of this Lemma depends on a $d$-dimensional version or if one wants, repeated application,
of the following Theorem.

\begin{theorem}[Mellin inversion theorem]\label{mellininversion}
Let $\phi(z)$ be function on $\bbbc$ satisfying the follwoing properties
\begin{itemize}
\item[(a)] $\phi$ is analytic in a vertical strip of the form $\alpha<x={\rm Re}(z)<\beta$
where $\alpha,\beta\in\bbbr$.
\item[(b)] $\int_{-\infty}^{\infty}|\phi(x+iy)|dy$ converges for all $x\in(\alpha,\beta)$.
\item[(c)] $\phi(z)\to0$ uniformly as $|y|\to\infty$ in $\alpha+\epsilon<x<\beta-\epsilon$
for all $\epsilon>0$.
\end{itemize}
Denote for all $t>0$,
$$f(t)={1\over2\pi i}\int_{c-i\infty}^{c+i\infty}t^{-z}\phi(z)dz.$$
Then,
$$\phi(z)=\int_0^{\infty}t^{z-1}f(t)dt.$$
\end{theorem}
For a proof of this theorem see \cite[Appendix 4, p341-342]{maclachlan}.
A treatment of the multidimensional case can also be found for example in
\cite{antipova}.
\medskip

\begin{proof} of Proposition \ref{mellinindependent}. Suppose we have a non-trivial relation $\sum_{\boldtheta\in\Theta}\lambda_{\boldtheta}M_{\boldtheta}=0$.
Let us use the notation $(\boldtheta\cdot B)(\v s)=\theta_1\v b_1\cdot\v s+\cdots
+\theta_N\v b_N\cdot\v s$. The
relation can be written as
$$0=\int_{\v s\in\boldsigma+\sqrt{-1}\bbbr^d}\left(\sum_{\boldtheta\in\Theta}
\lambda_{\boldtheta}e^{(\boldtheta\cdot B)(\v s)}\right)
|v_1|^{\v b_1\cdot\v s}\cdots|v_N|^{\v b_N\cdot\v s}\prod_{i=1}^N\Gamma(
-\gamma_i-\v b_i\cdot\v s)d\v s
$$
Let us now write $x_j=|v_1|^{b_{1j}}\cdots|v_N|^{b_{Nj}}$ where
$b_{ij}$ are the entries of the B-matrix. Then
$$|v_1|^{\v b_1\cdot\v s}\cdots|v_N|^{\v b_N\cdot\v s}=x_1^{s_1}\cdots x_d^{s_d}$$
By repeated use of the Mellin inversion Theorem \ref{mellininversion} we conclude
that the vanishing of the integral implies the identical vanishing of
$$\sum_{\boldtheta\in\Theta}
\lambda_{\boldtheta}e^{(\boldtheta\cdot B)(\v s)}$$
Since the exponentials are all distinct linear forms in $s_1,\ldots,s_d$
this implies that $\lambda_{\boldtheta}=0$ for all $\boldtheta\in\Theta$.
\end{proof}

Another proof of Proposition \ref{mellinindependent} can be found in \cite[Lemma 5.5]{nilsson}.
However, I hesitate somewhat about its completeness and decided to give the proof above.
\medskip

Since we assume that $\gamma_i<0$ for all $i$, all Mellin-Barnes integrals are solutions of
the corresponding A-hypergeometric system.
It would be very convenient if such a basis of solutions given by Mellin-Barnes integrals
would always exist. It turns out that with the exception of Appell $F_4$ all
$d=2$ systems Appell $F_1,F_2,F_3$  and Horn $G_1,G_2,G_3,H_1,\ldots,H_7$ this is the case.
A theoretical framework for a result like this may be provided in the
PhD-thesis of Lisa Nilsson \cite{nilsson}
suggesting that there does indeed exist such a basis if the complement of the so-called coamoeba of
the A-resultant is non-empty. In \cite{nilsson} this is elaborated for the case $d=2$.

Let us now make the following assumption on our system $H_A(\boldalpha)$.

\begin{assumption}\label{mellinbasis}
There exists a point $\v v_0\in(\bbbc^*)^N$ with an open neighbourhood in which
there exists a Mellin-Barnes basis of solutions.
\end{assumption}

For the practical determination of a Mellin-Barnes basis we use the following Proposition.

\begin{proposition}\label{anglechoice}
Let $H_A(\boldalpha)$ be a non-resonant system of rank $D$. The system allows a
Mellin-Barnes basis of solutions if and only if the zonotope $Z_B$ contains $D$ distinct
points $\boldtau_1,\ldots,\boldtau_D$ whose coordinates differ by integers.
Note that $Z_B$ is an open set in $\bbbr^d$.
\end{proposition}

\begin{proof}: From the discussion above it follows that the existence of a Mellin-Barnes basis
corresponds to the choice of $D$ $N$-tuples $(\theta_1,\ldots,\theta_N)$, representing argument
choices of a given point. Hence the differences between these $N$-tuples have coordinates which
are integer multiples of $2\pi$. The sums ${1\over2\pi}\sum_{i=1}^N\theta_i\v b_i$ are distinct,
hence the $D$ Mellin-Barnes basis elements correspond to $D$ points $\boldtau_i\in Z_B$
whose coordinates also differ by integers.

Suppose conversely we have $D$ points $\boldtau_i\in Z_B$ whose coordinates differ by integers.
Since the $\bbbz$-span of the columns $\v b_i$ is $\bbbz^d$ we can find for every $i$
integers $n_{i1},\ldots,n_{iN}$ such that $\tau_i-\tau_1=n_{i1}\v b_1+\cdots+n_{iN}\v b_N$.
So if $(\theta_1,\ldots,\theta_N)$ is an argument choice for $\tau_1$, then the $N$-tuples
$(\theta_1+2\pi n_{i1},\ldots,\theta_N+2\pi n_{iN})$ represent argument choices for $\tau_i$
with $i=1,2,\ldots,D$.
\end{proof}

For later use, we consider the vector $\boldtheta_j$ of arguments and
the vectors $\boldtau_j$ as column vectors. Then $2\pi\boldtau_j=B\boldtheta_j$ for
$j=1,2,\ldots,D$.

\section{Monodromy computation}\label{compute}
Suppose our system $H_A(\boldalpha)$ allows a Mellin-Barnes basis of solutions (Assumption \ref{mellinbasis}).
Denote the basis elements by $M_1,\ldots,M_D$ and the corresponding points in $Z_B$ by
$\boldtau_1,\ldots,\boldtau_D$. It is the goal of this section to compute the local monodromy groups with
respect to the basis $M_1,\ldots,M_D$. To this end we shall determine the transition matrices from
$M_1,\ldots,M_D$ to each of the bases of local series expansions.

To each point $\boldtau_i$ there corresponds a (not necessarily unique) choice
of arguments $\boldtheta_i=(\Arg(v_1),\ldots,\Arg(v_N))$ for $i=1,\ldots,D$.
We assume that the arguments are chosen such that
the differences $\boldtheta_i-\boldtheta_1$ have all their components equal to integer multiples of $2\pi$
(see Proposition \ref{anglechoice}).
Let $\v v_0\in(\bbbc^*)^N$ be a point whose coordinates have arguments $\boldtau_1$.
In particular we have a basis of Mellin-Barnes solutions around $\v v_0$. The Mellin-Barnes integral corresponding to the argument vector $\boldtheta_i$ is denoted by $M_i$.

Let $f_1,\ldots,f_D$ be a basis of series expansions for some local basis.
Suppose that the realm of convergence of these
local series expansions contains the torus $R$: $|v_i|=r_i$ for $i=1,\ldots,N$.
Choose a point $\v v_0'\in R$ with the same argument values as $\v v_0$
and let $\delta$ be a path from $\v v_0$ to $\v v_0'$ while keeping the
arguments fixed.
In a neighbourhood of $\v v_0'$ we also have a Mellin-Barnes basis of solutions
which are simply the analytic continuation of $M_1,\ldots,M_D$ along $\delta$.
For any $N$-tuple of integers $\v n=(n_1,\ldots,n_N)$ we consider the loop
$$c=c(\v n)=c(n_1,\ldots,n_N):
\ (e^{2\pi in_1t}v^{(0)'}_1,\ldots,e^{2\pi in_Nt}v^{(0)'}_N),
\ t\in[0,1].$$
Note that after analytic continuation of $M_1$ along the path $c((\boldtheta_j-\boldtheta_1)/2\pi)$ we
end up with the Mellin-Barnes solution $M_j$ for every $j$. We denote this path by $c_j$.

Let us denote a basis of local series expansions by $f_1,\ldots,f_D$ and the corresponding choices
of $\boldgamma$ by $\boldgamma^{(1)},\ldots,\boldgamma^{(D)}$. We regard the latter as row vectors.
Then there exist scalars $\mu_i$ such that
$$M_1=\mu_1f_1+\cdots+\mu_D f_D$$
in a neighbourhood of $\v v_0'$.
After continuation along $c_j$ the integral $M_1$ changes into $M_j$ for every $j$.
Under the paths $c_j$ the local expansions $f_i$ are multiplied by scalars. The space spanned by
$M_1,\ldots,M_D$ is $D$-dimensional. The space spanned by the images of $\mu_1 f_1+\cdots+\mu_D f_D$
under $c_j(\v v_0'),j=1,2,\ldots,D$ is at most equal to the number of non-zero $\mu_i$. Hence we conclude that
$\mu_i\ne0$ for all $i$. Let us renormalise the $f_i$ such that
$$M_1=f_1+\cdots+f_D.$$
Then after continuation along $c_j$ we get
$$M_j=e^{ i\boldgamma^{(1)}(\boldtheta_j-\boldtheta_1)}f_1+\cdots+e^{ i\boldgamma^{(D)}
(\boldtheta_j-\boldtheta_1)}f_D$$
for every $j$. Define
$$X_{\rho}=\pmatrix{1  & \cdots & 1\cr
e^{ i\boldgamma^{(1)}(\boldtheta_2-\boldtheta_1)} & \cdots &
e^{ i\boldgamma^{(D)}(\boldtheta_2-\boldtheta_1)}\cr
\vdots & & \vdots\cr
e^{ i\boldgamma^{(1)}(\boldtheta_D-\boldtheta_1)} & \cdots & e^{ i
\boldgamma^{(D)}(\boldtheta_D-\boldtheta_1)}\cr}.$$
Then
$$\pmatrix{M_1\cr \vdots\cr M_D\cr}=X_{\rho}\pmatrix{f_1\cr \vdots\cr f_D\cr}$$
hence $X_{\rho}$ is the desired transition matrix.
Let us now consider any closed path of the form $\delta^{-1}c(\v n)\delta,\v n\in\bbbz^N$ beginning and
ending in $\v v_0$. Continuation of $M_1,\ldots,M_D$ along $\delta$ is trivial since
the Mellin-Barnes integrals converge throughout. However these integrals do not converge anymore
if we continue along $c(\v n)$. For that we have to change to the local basis $f_i,i=1,\ldots,D$.
Analytic continuation along $c(\v n)$ changes them into $e^{2\pi i\v n\cdot\boldgamma^{(i)}}f_i$ for $i=1,\ldots,D$.
Express these solutions in terms of the $M_i$ again and continue back along $\delta^{-1}$.
The monodromy matrix can be computed as follows.
Let $\chi_{\rho}(\v n)$ be the diagonal $D\times D$-matrix
with entries $e^{2\pi i\boldgamma^{(i)}\v n},\ i=1,\ldots,D$. It is the monodromy matrix with respect to $f_1,\ldots,f_D$.
With respect to $M_1,\ldots,M_D$ this monodromy element has matrix $X_{\rho}\chi_{\rho}(\v n) X_{\rho}^{-1}$.
Thus we see that all local monodromies can be written with respect to a fixed Mellin-Barnes basis.

\section{An implementation}\label{algorithm}
The considerations in the previous sections, together with some practical tricks,
lead to an algorithm to compute monodromy matrices, which we describe in this
section.

We start with a totally non-resonant hypergeometric system $H_A(\boldalpha)$ and we
assume that there exists a Mellin-Barnes basis.
The starting data are a $d\times N$ B-matrix $B$ and a parametervector $\boldgamma_0$
(in row form) such that $A\boldgamma_0=\boldalpha$.
In general both $B$ and $\boldgamma_0$ can easily be read off
from an explicit series solution of a hypergeometric system.
For example, from the expansion of Appell $F_2$ as on page \pageref{F2expansion2}.
We also assume we know the rank $D$ of the system.

\smallskip

{\it Step} 1. Using the B-matrix we determine the zonotope $Z_B$ and a find
$D$ distinct points in it, whose coordinates differ by integers. Since we assumed
the existence of a Mellin-Barnes basis these points exist.
Call the points $\boldtau_1,\ldots,\boldtau_D$. From the proof of Proposition
\ref{anglechoice} we know that to each $\tau_i$ there exists a column vector of arguments
$\theta_i\in\bbbr^n$ such that $2\pi\tau_i=B\theta_i$. However, we do not compute
these angle vectors.
\smallskip

{\it Step} 2. Construct the set ${\cal I}$ of all subsets $I$ of cardinality $d$
of the columns $\{\v b_1,\ldots,\v b_N\}$
with $\delta_I=|\det_{i\in I}(\v b_i)|\ne0$. As a fine point, if $\Delta_I>1$ we include
$\Delta_I$ copies of $I$ in ${\cal I}$. For each $I$ there exists a parametervector
$\boldgamma^I$ in the following way.
Denote the rows of the B-matrix $B$ by $\v l_1,\ldots,\v l_d$,
recall that this is a basis of the lattice $L$.
In case $\Delta_I=1$ we take the uniquely determined real
numbers $\mu_1,\ldots,\mu_d$ such that $\boldgamma_0+\mu_1\v l_1+\cdots+\mu_d\v l_d$
has $i$-th coordinate 0 for all $i\in I$ and call this sum $\gamma^I$.
In case $\Delta_I>1$ we make $\Delta_I$ choices for $(\mu_1,\ldots,\mu_d)$, distinct
modulo $\bbbz^d$, such that $\boldgamma_0+\mu_1\v l_1+\cdots+\mu_d\v l_d$ has integer
coordinates on the $i$-th position for all $i\in I$. In this way we get $\Delta_I$
different parametervectors $\boldgamma^I$ (in row form) for a given $I\in{\cal I}$.
However, in the computation we only retain the (row vectors) $\boldmu^I=(\mu_1,\ldots,\mu_d)^t$
for each $I$.
They have the property that $\boldgamma^I-\boldgamma_0=\boldmu^I B$ for all $I$.
\smallskip

{\it Step} 3. To every $I$ we associate a column vector $X_I$ of length $D$
given by
\begin{equation}\label{Xvector}
X_I=(1,\exp(2\pi i\boldmu^I(\boldtau_2-\boldtau_1)),
\ldots,\exp(2\pi i\boldmu^I(\boldtau_D-\boldtau_1)))^t.
\end{equation}
Since, for each $j=1,2,\ldots,D$ we have
\begin{eqnarray*}
2\pi \boldmu^I(\boldtau_j-\boldtau_1)&=&\mu^IB(\theta_j-\theta_1)\\
&=&(\boldgamma^I-\boldgamma_0)(\theta_j-\theta_1)
\end{eqnarray*}
the $j$-th components of $X_I$ differs by a factor
$$\exp(-i\boldgamma_0(\theta_j-\theta_1))$$
from the similar components in the columns of the transition matrices $X_{\boldrho}$.
The only effect is that the transition matrices built out of our present $X_I$
will give us the transition of the Mellin-Barnes basis to a renormalized local
basis. This will have no effect on the monodromy computation.
\smallskip

{\it Step} 4. For every cone in the secondary fan (specified by a convergence
direction $\rho$ inside that cone) we determine the sets $I\in{\cal I}$ such that
the cone or, equivalently, $\rho$ lies in the positive real cone spanned by the
vectors $\{\v b_i\}_{i\in I}$. Call this set of sets ${\cal I}_{\rho}$.
The theory of Gel'fand, Kapranov and Zelevinsky tells
us that ${\cal I}_{\rho}$ contains precisely $D$ sets (when we count possible repetitions
of a set with $\Delta_I>1$). Let $X_{\rho}$ be the $D\times D$-matrix
whose columns are the vectors $X_I$ with $I\in{\cal I}_{\rho}$.
The matrices $X_{\rho}$ are the
transition matrices from the Mellin-Barnes basis to the local power series basis,
all of whose elements contain $\rho$ as a convergence direction.
\smallskip

{\it Step} 5. For every cone in the secondary fan (specified by a convergence
direction $\rho$) we determine the characters of the local monodromies of the corresponding
power series solutions. Let $I_0$ be such that $|\det(\v b_j)_{j\in I_0}|=1$
We choose a set of $d$ generators
for the local monodromy as follows. For every $j=1,\ldots,d$ we define $\v n_j\in\bbbz^N$
such that $\v n_j$ has support in $I_0$ and $B\v n_j=\v e_j$, the $j$-th standard basis
vector in $\bbbr^d$. Since the support is in $I_0$ we have $\boldgamma_0\v n_j=0$.
Furthermore,
$$\boldgamma^I\v n_j=(\boldgamma^I-\boldgamma_0)\v n_j=
(\boldmu^I-\boldmu^{I_0})B\v n_j=\mu^I_j-\mu^{I_0}_j$$
for $j=1,\ldots,d$. Hence the characters corresponding to the path $c(\v n_j)$ read
$$\exp(2\pi i\boldgamma^I\v n_j)=\exp(2\pi i(\mu^I_j-\mu^{I_0}_j)).$$
\smallskip

{\it Step} 6. This is the final step in which we compute $d$ monodromy matrices
for every cone of the secondary fan. For a cone, specified by a convergence direction
$\rho$, and a loop $c(\v n_j)$ (defined in step 5) we construct the matrix $X_{\rho}$ as in Step 4, and a
diagonal matrix $\chi_{\boldrho,j}$ with entries $\exp(2\pi i\mu^I_j),\ I\in{\cal I}_{\rho}$,
as in Step 5.
We see to it that both in $X_{\rho}$ and $\chi_{\boldrho,j}$ we keep the same ordering
of the set ${\cal I}_{\rho}$. Then construct the matrix
$$M_{\boldrho,j}=X_{\rho}\chi_{\rho,j}X_{\rho}^{-1}.$$
\smallskip

Let $F$ be the number of open cones in the secondary fan. Then we get $dF$ monodromy
matrices in this way. They generate a subgroup of the monodromy group whose
projectivization (quotient by scalars) we denote by $lMon$. As remarked before,
computation of local monodromies by other loops will
only add scalar matrices and therefore does not change $lMon$.

\section{An example, Appell $F_2$}
Recall that a B-matrix is given by
$$B=\pmatrix{-1 & -1 & 0 & 1 & 0 & 1 & 0\cr -1 & 0 & -1 & 0 & 1 & 0 & 1\cr}$$
and a parameter vector $\boldgamma_0=(-\alpha,-\beta,-\beta',\gamma-1,\gamma'-1,0,0)$.
We trust that no
confusion will arise with the existing notations $\boldalpha$ and $\boldgamma$.
Two powerseries solution expansions have already been given on pages \pageref{F2expansion1}
and \pageref{F2expansion2}.
The set ${\cal I}$ consists of 15 elements, $\{1,3\}$, $\{1,4\}$, $\{1,6\}$,
$\{1,5\}$, $\{1,7\}$, $\{1,2\}$, $\{2,3\}$, $\{2,5\}$, $\{2,7\}$, $\{3,4\}$, $\{3,6\}$,
$\{4,5\}$, $\{4,7\}$, $\{5,6\}$, $\{6,7\}$. Here is a table with the corresponding values of $\mu^J$.

\begin{center}
\begin{tabular}{|l|c|l|}
\hline
&$J$ & $\boldmu^J$\\
\hline
1 & $\{1,3\}$ & $ -\alpha + \beta', -\beta'$ \\
2 & $\{1,4\}$ & $1 - \gamma, -1 - \alpha + \gamma$  \\
3 & $\{1,6\}$ & $0, -\alpha$  \\
4 & $\{1,5\}$ & $-1 - \alpha + \gamma', 1 - \gamma'$\\
5 & $\{1,7\}$ & $-1 + \gamma', -\alpha, 0$ \\
6 & $\{1,2\}$ & $-\beta, -\alpha + \beta$  \\
7 & $\{2,3\}$ & $-\beta, -\beta'$\\
\hline
\end{tabular}$\qquad$
\begin{tabular}{|l|c|l|}
\hline
8 & $\{2,5\}$ & $-\beta, 1 - \gamma'$  \\
9 & $\{2,7\}$ & $-\beta, 0$\\
10 & $\{3,4\}$ & $1 - \gamma, -\beta'$  \\
11 & $\{3,6\}$ & $0, -\beta'$\\
12 & $\{4,5\}$ & $1 - \gamma, 1 - \gamma'$\\
13 & $\{4,7\}$ & $1 - \gamma, 0$\\
14 & $\{5,6\}$ & $0, 1 - \gamma'$\\
15 & $\{6,7\}$ & $0, 0$\\
\hline
\end{tabular}
\end{center}

As noted earlier, the zonotope $Z_B$ contains the four points
$(\pm1/2,\pm1/2)$. Define
$$\boldtau_1=(-1/2,-1/2)^t,\ \boldtau_2=(1/2,-1/2)^t,\ \boldtau_3=(-1/2,1/2)^t,
\ \boldtau_4=(1/2,1/2)^t.$$
Then the vectors $X_J$, as defined in Step 3 of our algorithm on page \pageref{Xvector},
read $$1,e(\mu^J_1),e(\mu^J_2),e(\mu^J_1+\mu^J_2)$$
where we use the notations $e(x)=e^{2\pi ix}$ and
$a=e(\alpha),b=e(\beta),b'=e(\beta'),c=e(\gamma),c'=e(\gamma')$.
Here is the list of all $X_J$ with the same ordering as in the previous table,
\begin{center}\label{listofXvectors}
\begin{tabular}{|l|l|cccc|}
\hline
&$J$ & & &$X_J$\hfill& \\
\hline
1 & $\{1,3\}$ & $1$ & $b'/a$ & $1/b'$   & $1/a$  \\
2 & $\{1,4\}$ & $1$ & $1/c$  & $c/a$    & $1/a $ \\
3 & $\{1,6\}$ & $1$ & $1$    & $1/a$    & $1/a$ \\
4 & $\{1,5\}$ & $1$ & $c'/a$ & $1/c')$  & $1/a$\\
5 & $\{1,7\}$ & $1$ & $1/a$  & $1$      & $1/a$\\
6 & $\{1,2\}$ & $1$ & $1/b$   & $b/a$   & $1/a$ \\
7 & $\{2,3\}$ & $1$ & $1/b$   & $1/b'$   & $1/bb'$\\
\hline
\end{tabular}$\qquad$
\begin{tabular}{|l|l|cccc|}
\hline
8 & $\{2,5\}$ & $1$ & $1/b$   & $1/c'$   & $1/bc'$ \\
9 & $\{2,7\}$ & $1$ & $1/b$   & $1$      & $1/b$\\
10 & $\{3,4\}$ & $1$ & $1/c$  & $1/b'$   & $1/cb'$ \\
11 & $\{3,6\}$ & $1$ & $1$    & $1/b'$   & $1/b'$ \\
12 & $\{4,5\}$ & $1$ & $1/c$  & $1/c'$   & $1/cc'$ \\
13 & $\{4,7\}$ & $1$ & $1/c$  & $1$      & $1/c$ \\
14 & $\{5,6\}$ & $1$ & $1$    & $1/c'$   & $1/c'$ \\
15 & $\{6,7\}$ & $1$ & $1$    & $1$      & $1$\\
\hline
\end{tabular}
\end{center}

To write down local monodromies we use Step 5 of our algorithm. The characters
$e(\mu_1)$ are said to correspond to path I and the characters $e(\mu_2)$ correspond to path II.
We do not need to write down a separate table for them since they are simply the
second and third component of the vectors $X_J$.

As an example for the action of path I on a local basis we take the four local basis
solutions with convergence direction $-0.5,1$, as before.
The transition matrix $X_{\rho}$ consists of the vectors $X_J$ with numbering
4,5,8,9 in Table \ref{listofXvectors}.
The transition matrix reads
$$
X_{\rho}=\pmatrix{1 & 1 & 1 & 1\cr  c'/a & 1/a & 1/b & 1/b\cr
 1/c' & 1 & 1/c' & 1 \cr 1/a & 1/a & 1/bc' & 1/b\cr}
$$
For the path I we get the monodromy matrix
$$X_{\rho}\pmatrix{c'/a & 0 & 0 & 0\cr 0 & 1/a & 0 & 0\cr
0 & 0 & 1/b & 0\cr 0 & 0 & 0 & 1/b\cr}X_{\rho}^{-1}$$
which equals
$$\pmatrix{0 & 1 & 0 & 0\cr (-1+c')/ab & 1/b+1/a+c'/a &
c'/ab & -c'/a\cr 0 & 0 & 0 & 1\cr
-1/ab & 1/a & 0 & 1/b\cr}
$$
with respect to the Mellin-Barnes basis. For the path II we get
$$\pmatrix{0 & 0 & 1 & 0\cr 0 & 0 & 0 & 1\cr -1/c' & 0 & 1+1/c' & 0\cr
0 &-1/c' & 0 & 1+1/c'\cr}.$$

The calculation so far has been carried out for the convergence direction
$(-0.5,1)$. In fact we get the same matrices for every convergence direction in
the cone spanned by $\v b_2,\v b_5$ of the secondary fan.
We can proceed in the same way with the other four cones. In each case we find
two monodromy matrices. After removing duplicate matrices we end up with six
monodromy matrices. They are given in the Appendix of this paper, together with
a comparison of them with the five generators given by M.Kato in \cite{kato}.
It turns out that the group generated by our six generators is conjugate
to the group computed in \cite{kato}.

\section{An example, Clausen $_3F_2$}\label{example3F2}
In this section we apply our method to the case of one variable
$$\hypergeo{3}{2}{\alpha_1,\alpha_2,\alpha_3}{\beta_1,\beta_2}{z}$$
which, up to a constant factor, is defined by the series
$$\sum_{n\ge0}{\Gamma(\alpha_1+n)\Gamma(\alpha_2+n)\Gamma(\alpha_3+n)
\over\Gamma(\beta_1+n)\Gamma(\beta_2+n)n!}\ z^n.$$
Using the identity $\Gamma(z)\Gamma(1-z)=\pi/\sin\pi z$ we see that the
series is proportional to
$$\sum_n{(-z)^n\over\Gamma(1-\alpha_1-n)\Gamma(1-\alpha_2-n)\Gamma(1-\alpha_3-n)
\Gamma(\beta_1+n)\Gamma(\beta_2+n)\Gamma(1+n)}.$$
So the B-matrix is given by
$$B=(-1,-1,-1,1,1,1)^t$$
and $Z_B$ is simply the open interval $(-3/2,3/2)$. In it we can take the three
points $\tau_1=-1,\tau_2=0,\tau_3=1$ and so we see that we have a Mellin-Barnes basis of solutions.
For the set $I_0$ we take $\{6\}$ and
$$\boldgamma_0=(-\alpha_1,-\alpha_2,\alpha_3,\beta_1,\beta_2,0).$$
We consider the components modulo $\bbbz$.
The set of columns of $B$ has 6 subsets of cardinality 1 and the corresponding values
of $\mu_1$ are
$$\alpha_1,\alpha_2,\alpha_3,-\beta_1,-\beta_2,0.$$
Letting $a_i=e(\alpha_i)$ and $b_j=e(\beta_j)$ we get for the vectors $X_J$,
$$(1,a_1,a_1^2),\ (1,a_2,a_2^2),\ (1,a_3,a_3^2)$$
$$(1,b_1,b_1^2),\ (1,b_2,b_2^2),\ (1,1,1).$$
There is only one loop to consider for every local basis. Consider the convergence direction $-1$.
This lies in the positive cones spanned by $\v b_1,\v b_2,\v b_3$ respectively.
The transition matrix reads
$$X_{\rho}=\pmatrix{1 & 1 & 1\cr a_1 & a_2 & a_3\cr a_1^2 & a_2^2 & a_3^2\cr}$$
and the diagonal character matrix
$$\chi_{\rho}=\pmatrix{a_1 & 0 & 0\cr 0 & a_2 & 0\cr 0 & 0 & a_3\cr}.$$
We get
$$X_{\rho}\chi_{\rho}X_{\rho}^{-1}=
\pmatrix{0 & 1 & 0\cr 0 & 0 & 1\cr a_1a_2a_3 & -a_2a_3-a_2a_1-a_3a_1 & a_1+a_2+a_3}.$$
This is precisely the matrix representation for the monodromy matrix around $z=\infty$
for $_3F_2$ as given in \cite{beukersheckman}. We get a similar result for the
monodromy matrix around $z=0$ (with $b_1,b_2,1$ instead of $a_1,a_2,a_3$)

\section{Existence of Mellin-Barnes bases}

In this section we show that certain families of hypergeometric equations satify
Assumption \ref{mellinbasis}, and some don't (the case of Lauricell $F_C$).

\subsection{Lauricella $F_A$}
The Lauricella system $F_A$ in $n$ variables is a system of rank $2^n$.
From the powerseries
$$F_A(a,\v b,\v c|\v x)=\sum_{\v m\ge 0}
	{(a)_{|\v m|}(\v b)_{\v m}\over (\v c)_{\v m}\v m!}
	\v x^{\v m}$$
in $\v x=(x_1,\ldots,x_n)$ we see that an $n\times(3n+2)$ B-matrix is given by
$$\pmatrix{
1 & 1 & 0 &\cdots & 0 & -1 & 0 & \cdots & 0 & -1 & 0 & \cdots 0\cr
1 & 0 & 1 &\cdots & 0 & 0 & -1 & \cdots & 0 & 0 & -1 & \cdots 0\cr
\vdots&&&&&&&&&&&\vdots\cr
1 & 0 & 0 &\cdots & 1 & 0 & 0 & \cdots & 0 & 0 & 0 & \cdots -1\cr}
$$
The B-zonotope is thus given by the points
$$\lambda (\v e_1+\cdots\v e_n)+\sum_{i=1}^m\mu_i\v e_i$$
where $|\lambda|<1/4$ and $|\mu_i|<3/4$ for $i=1,\ldots,n$.
Let us choose $\epsilon>0$ sufficiently small. Consider the
$2^n$ points
$$(3/4-2\epsilon)(\v e_1+\cdots+\v e_n)-k_1\v e_1-\cdots k_n\v e_n$$
where $k_i\in\{0,1\}$ for all $i$. Each such point equals
$$(1/4-\epsilon)(\v e_1+\cdots+\v e_n)+(1/2-k_1-\epsilon)\v e_1+
\cdots+(1/2-k_n-\epsilon)\v e_n$$
which is clearly contained in the B-zonotope.
\medskip

\subsection{Lauricella $F_B$}
The Lauricella system $F_B$ is also a system of rank $2^n$.
From the powerseries
$$F_B(\v a,\v b, c|\v x)=\sum_{\v m\ge 0}
	{(\v a)_{\v m}(\v b)_{\v m}\over (c)_{|\v m|}\v m!}
	\v x^{\v m}$$
in $\v x=(x_1,\ldots,x_n)$ we see that an $n\times(3n+2)$ B-matrix is given by
$$\pmatrix{
1 & 0 &\cdots & 0 & 1 & 0 & \cdots & 0 & -1 & -1 & 0 & \cdots 0\cr
0 & 1 &\cdots & 0 & 0 & 1 & \cdots & 0 & -1 & 0 & -1 & \cdots 0\cr
\vdots&&&&&&&&&&&\vdots\cr
0 & 0 &\cdots & 1 & 0 & 0 & \cdots & 1 & -1 & 0 & 0 & \cdots -1\cr}
$$
Note the B-zonotope is the same as in the case of Lauricella $F_A$. Hence
there exists a Mellin-Barnes basis of solutions.
\medskip

\subsection{Lauricella $F_D$}
The Lauricella system $F_D$ in $n$ variables is a system of rank $n+1$.
From the powerseries
$$F_D(a,\v b, c|\v x)=\sum_{\v m\ge 0}
	{(a)_{|\v m|}(\v b)_{\v m}\over (c)_{|\v m|}\v m!}
	\v x^{\v m}$$
in $\v x=(x_1,\ldots,x_n)$ we deduce an $n\times(2n+2)$ B-matrix
$$\pmatrix{
1 & 1 & 0 & \cdots & 0 & -1 & -1 & 0 & \cdots & 0\cr
1 & 0 & 1 & \cdots & 0 & -1 & 0 & -1 & \cdots & 0\cr
\vdots &&&&&&&&&\vdots\cr
1 & 0 & 0 & \cdots & 1 & -1 & 0 & 0 & \cdots & -1\cr
}$$
Hence the B-zonotope consists of the points
$$\lambda_0(\v e_1+\cdots+\v e_n)+\lambda_1\v e_1+\cdots+\lambda_n\v e_n$$
where $|\lambda_i|<1/2$ for $i=0,1,\ldots,n$.
Choose $\epsilon>0$ sufficiently small and consider the $n+1$ points
$$-\epsilon(n\v e_1+(n-1)\v e_2+\cdots+2\v e_{n-1}+\v e_n)+\sum_{i=0}^k\v e_i$$
for $k=0,1,2\ldots,n$. Each such point can be rewritten as
$$(1/2-(n-k-1/2)\epsilon)(\v e_1+\cdots+\v e_n)+
\sum_{j=1}^n(\pm1/2-(k-j-1/2)\epsilon)\v e_j$$
where $\pm1/2$ is $1/2$ if $k>j+1/2$ and $-1/2$ if $k< j+1/2$.
Hence they are contained in the B-zonotope and we have found a
Mellin-Barnes basis for Lauricella $F_D$.
\medskip

\subsection{Lauricella $F_C$}
The Lauricella system $F_C$ in $n$ variables is system of rank $2^n$.
From the powerseries
$$F_C(a, b, \v c|\v x)=\sum_{\v m\ge 0}
	{(a)_{|\v m|}(\v b)_{|\v m|}\over (c)_{\v m}\v m!}
	\v x^{\v m}$$
in $\v x=(x_1,\ldots,x_n$ we deduce an $n\times(2n+2)$ B-matrix
$$\pmatrix{
1 & 1 & -1 & 0 & \cdots & 0 & -1 & 0 & \cdots & 0\cr
1 & 1 & 0 & -1 & \cdots & 0 & 0 & -1 & \cdots & 0\cr
\vdots &&&&&&&&&\vdots\cr
1 & 1 & 0 & 0 & \cdots & -1 & 0 & 0 & \cdots & -1\cr
}$$
Note that the B-zonotope is the same as for Lauricella $F_D$, but
this time we have to find a Mellin-Barnes basis of $2^n$ solutions.
Clearly this is impossible if $n>1$.

\subsection{Aomoto-Gel'fand system $E(3,6)$}
This system forms the subject of the second part of
M.Yoshida's book \cite{yoshida}. It is an Aomoto system which can be reinterpreted as
an A-hypergeometric system. It has four essential variables ($d=4$) and rank 6.
 system that corresponds to configurations
of six points (or lines) in $\bbbp^2$. We start by giving the A-matrix of the system,
$$A=\pmatrix{
0 & 1 & 0 & 0 & 1 & 0 & 0 & 1 & 0 \cr
0 & 0 & 1 & 0 & 0 & 1 & 0 & 0 & 1\cr
1 & 1 & 1 & 0 & 0 & 0 & 0 & 0 & 0\cr
0 & 0 & 0 & 1 & 1 & 1 & 0 & 0 & 0\cr
0 & 0 & 0 & 0 & 0 & 0 & 1 & 1 & 1\cr}$$
and the parameters $(\alpha_1,\alpha_2,2-\alpha_4,2-\alpha_5,2-\alpha_6)$.
The integrand of the Euler integral as defined in \cite[p 607]{beukersalgebraic}
reads
$${s^{\alpha_1-1}t^{\alpha_2-1}t_1^{1-\alpha_4}t_2^{1-\alpha_5}t_3^{1-\alpha_6}
\over 1-t_1(s+t+1)-t_2(s+v_1t+v_3)-t_3(s+v_2t+v_4)}
ds\wedge dt\wedge dt_1\wedge dt_2\wedge dt_3.$$
Perform the substitutions $t_1\to t_1/(s+t+1), t_2\to t_2/(s+v_1t+v_3)$ and
$t_3\to t_3/(s+v_2t+v_4)$. We obtain the integrand
$${t_1^{1-\alpha_4}t_2^{1-\alpha_5}t_3^{1-\alpha_6}
\over 1-t_1-t_2-t_3}\prod_{i=1}^6(L_i)^{\alpha_i-1}
ds\wedge dt\wedge dt_1\wedge dt_2\wedge dt_3$$
where
$$L_1=s,\ L_2=t,\ L_3=1,\ L_4=s+t+1,$$
$$L_5=s+v_1t+v_3,\ L_6=s+v_2t+v_4$$
and $\alpha_3=3-\alpha_1-\alpha_2-\alpha_4-\alpha_5-\alpha_6$.
Integration with respect to $t_1,t_2,t_3$ leaves us with a 2-form which is
the integrand given on page 221 of \cite{yoshida}, but with $v_i$ instead of $x^i$.
Thus we see that our A-matrix corresponds to a Gel'fand-Aomoto system which is
associated to configurations of six lines in $\bbbp^2$. The system is a
four variable system of rank 6. It is irreducible if and only if none of
the $\alpha_i$ is an integer, see \cite[Prop 2]{masatayo2}. A possible B-matrix reads
$$B=\pmatrix{
1& 0& -1& 0& 0& 0& -1& 0& 1\cr
1& -1& 0& 0& 0& 0& -1& 1& 0\cr
1& 0& -1& -1& 0& 1& 0& 0& 0\cr
1& -1& 0& -1& 1& 0& 0& 0& 0\cr}^t$$
The set ${\cal I}$ consists of 81 sets, and hence 81 distinct local solutions. The
number of local solution bases is 108.
In a straightforward manner one can check that the B-zonotope $Z_B$ contains
the points
$$p,p+(0,0,0,1),p+(1,0,0,0),p+(1,0,1,1),p+(1,1,0,1),p+(1,1,1,1)$$
where $p=(-0.9,-0.4,-0.5,-0.7)$. Hence the system $E(3,6)$ has a Mellin-Barnes
basis of solutions. Computation of a set of generators for the monodromy group
$lMon$ is now straightforward. We get $82$ matrices, but have not made an attempt
to compare with the $20$ generators of $Mon$ found in \cite{masatayo1}.

\section{Hermitian forms}
In the cases where we carried out the algorithm given above, it turns out that whenever
$\boldalpha\in\bbbr^r$, and the system is totally nonresonant,
there exists a unique (up to a constant factor) hermitean form
which is invariant under the group $lMon$. Subsequent studies lead us to the following
conjecture.

\begin{conjecture}\label{invariantform}
Let $H_A(\boldalpha)$ be a non-resonant A-hypergeometric system with $\boldalpha\in\bbbr^r$.
Then there exists a non-trivial unique (up to scalars)
Hermitean form, invariant under the monodromy group. More concretely, there exist a
Hermitean $D\times D$-matrix $H$ such that $\overline{g}^tHg=H$ for all elements $g$ of
the monodromy group. Here $D$ denotes the rank of $H_A(\boldalpha)$.

Moreover, when the system is totally non-resonant,
the signature of $H$ is determined by the signs of the numbers
$$\prod_{i\not\in I}\sin(\pi\gamma_i^I)$$
as $I$ runs through the
elements of ${\cal I}_{\boldrho}$ for some convergence direction $\rho$.
\end{conjecture}

We want to deal with this matter in another paper. However, we do like to add
that a detailed calculation shows that the signatures thus obtained are in accordance
with the results on $E(3,6)$ in \cite[Prop 1]{masatayo2} (except for a small printing error).
Note that signature $(5,1)$ does not occur. Similarly, calculations for Lauricella $F_D$
give us results which are in accordance with Picard \cite{picard}, Terada \cite{terada}
and Deligne-Mostow \cite{delignemostow}.

\ifx
\begin{proof} (sketch):
In \cite{beukersalgebraic} we find that the solution space can be given by integrals
over suitable cycles of integration (generalised Pochhammer cycles) or over twisted (or loaded)
cycles in the sense of Kita and Yoshida \cite{kitayoshida}.
The Hermitean form arises from the intersection form on these cycles.
\end{proof}

In this section we shall assume that Remark \ref{invariantform holds}.

Suppose that $\boldalpha\in\bbbr^r$ and that $H_A(\boldalpha)$ is totally non-resonant. Suppose also that we have a Mellin-Barnes basis of solutions. Let us recall the notations
$\v b_I,{\cal I}$ from Section \ref{powerseries}. To each
$I\in{\cal I}$ there is associated a choice of $\boldgamma$ which we denote $\boldgamma^I$.
Let $\boldtheta_1,\ldots,\boldtheta_D\in \bbbr^N$ be the argument vectors corresponding
to the Mellin-Barnes basis elements $M_1,\ldots,M_D$. To each $I\in{\cal I}$ we associate the
vector
\begin{equation}
X_I=(1,\exp( i\ (\boldtheta_2-\boldtheta_1)\cdot\boldgamma^I),\ldots,\exp( i\ (\boldtheta_D-\boldtheta_1)\cdot\boldgamma^I)^t\in\bbbc^D,
\end{equation}
as on page \pageref{Xvector}. Vectors of this form are the column vectors of the
transition matrices $X_{\rho}$ from Section \ref{compute}.

To every monodromy element we associate $g\in GL(n,\bbbc)$ such that the monodromy substition
has the form $(M_1,\ldots,M_D)^t\to g(M_1,\ldots,M_D)^t$. By Theorem \ref{invariantform}
there exists a Hermitean $n\times n$-matrix
$H$ (i.e. $\overline{H}^t=H$) such that $\overline{g}^tHg=H$ for all monodromy matrices $g$.

\begin{proposition}\label{orthogonal}
Let notations be as above. Let $I,J\in{\cal I}$ and suppose that $\v b_I$ and $\v b_J$
have a point in common. Then $X_I,X_J$ are orthogonal with respect to $H$, i.e $\overline{X}_J^tHX_I=0$.
\end{proposition}

\begin{proof}: Given a basis of local solutions we have seen in
Section \ref{compute} that the local monodromies with respect
to the Mellin-Barnes basis have the form $g=X_{\rho}\chi X_{\rho}^{-1}$ where $\chi$ is a diagonal matrix representation
with entries
$$\v n\in\bbbz^N\mapsto (\exp(2\pi  i\ \v n\cdot\boldgamma^{(1)}),
\ldots,\exp(2\pi  i\ \v n\cdot\boldgamma^{(D)})).$$
We first show that all these characters are distinct. Suppose $\exp(2\pi i\ \v n\cdot\boldgamma^{(k)})=
\exp(2\pi i\ \v n\cdot\boldgamma^{(l)})$ for all $\v n\in\bbbz^N$. This implies that $\boldgamma^{(k)}\is\boldgamma^{(l)}
\mod{\bbbz^N}$. Since also $\boldgamma^{(k)}-\boldgamma^{(l)}\in L\otimes\bbbr$ we conclude that
$\boldgamma^{(k)}-\boldgamma^{(l)}\in L$ and hence $k=l$.

Let $\rho$ be a common point of $\v b_I$ and $\v b_J$ and consider the corresponding
basis of local power series solutions indexed by $I_{\rho}$.
With the notation from Section \ref{compute} we proceed by rewriting $\overline{g}^tHg=H$ as
$$(\overline{X}_{\rho}^t)^{-1}\chi^{-1}\overline{X}_{\rho}^tHX_{\rho}\chi X_{\rho}^{-1}=H.$$
This implies $\chi^{-1}\overline{X}_{\rho}^t HX_{\rho}\chi=\overline{X}_{\rho}^tHX_{\rho}$. Since $\chi$
consists of distinct characters on its diagonal we conclude that $\overline{X}_{\rho}^tHX_{\rho}$
is a diagonal matrix. Hence the columns of $X_{\rho}$ are orthogonal with respect to $H$.
In particular $X_I,X_J$ are orthogonal with respect to $H$.
\end{proof}

A question is whether the orthogonality property of Proposition \ref{orthogonal} suffices to determine $H$ uniquely (up to a scalar factor).
Also, to determine $H$ completely we need to know the values of
$\overline{X}^t_IHX_I$ for every $I$. So far we have not been able to answer these questions.
Instead we end with the following conjecture.

\begin{conjecture}\label{signature}
After suitable normalisation of $H$ we have
for every $I\in{\cal I}$ that $\overline{X}^t_IHX_I=\Delta_I\prod_{i\not\in I}\sin(\pi\gamma_i^I)$.
\end{conjecture}

Consequently, this conjecture implies that the {\it signature of the hermitian form} $H$ is determined
by the signs of the values of $\prod_{i\not\in I}\sin(\pi\gamma_i^I)$ as $I$ runs through the
elements of ${\cal I}_{\boldrho}$ for some convergence direction $\rho$.

\fi

\section{Appendix}
Here we reproduce the six matrices obtained from the monodromy calculation of
Appell $F_2$.
$$M_1=\pmatrix{0 & 0 & 1 & 0\cr 0& 0& 0 & 1\cr
-(1+c)/ab')& c/ab' & 1/b'+1/a+c/a&-c/a\cr
-1/ab' & 0 & 1/a & 1/b'\cr},$$

$$M_2=\pmatrix{0 & 0 & 1 & 0\cr 0 & 0 & 0 & 1\cr
-b/ab'-1/c' & b/c' & 1+b/a+1/c' & b(-1+1/b'-
1/c')\cr -1/ab' & 0 & 1/a & 1/b'\cr},$$

$$M_3=\pmatrix{0 & 0 & 1 & 0\cr 0 & 0 & 0 & 1\cr
-1/c' & 0 & 1+1/c' & 0\cr 0 & -1/c' & 0 & 1+1/c'\cr},$$

$$M_4=\pmatrix{0 & 1 & 0 & 0\cr
-b'/ab-1/c & 1+b'/a+1/c & b'/c &b'
(-1+1/b-1/c)\cr
0 & 0 & 0 & 1\cr -1/ab & 1/a & 0 & 1/b\cr},$$

$$M_5=\pmatrix{0 & 1 & 0 & 0\cr -1/c & 1+1/c & 0 & 0\cr
0 & 0 & 0 & 1\cr 0 & 0 & -1/c & 1+1/c\cr},$$

$$M_6=\pmatrix{0 & 1 & 0 & 0\cr
-(1+c')/ab & 1/b+1/a+c'/a & c'/ab & -c'/a\cr
0 & 0 & 0 & 1\cr -1/ab & 1/a & 0 & 1/b\cr}.$$

Now let $g_1,g_2,g_3,g_4,g_5$ be the monodromy matrices defined in formulas
(2.7), (2.8), (2.9), (2.10), (2.11) in M.Kato's paper \cite{kato}, where
our symbols $a,b,b',c,ç'$ are Kato's symbols $e(a),e(b),e(b'),e(c),e(c')$.
Define the conjugation matrix
$$S=\pmatrix{-1 & c & c' & -cc'\cr
-1 & 1 & c' & -c'\cr -1 & c & 1 & -c\cr -1 & 1 & 1 & -1\cr}.$$
Then the relations between the $M_i$ and $g_j$ are given by
$$M_1=S^{-1}g_2g_3g_5S\quad M_2=S^{-1}g_2g_5S,\quad M_3=S^{-1}g_2S$$
$$M_4=S^{-1}g_1g_4S,\quad M_5=S^{-1}g_1S,\quad M_6=S^{-1}g_1g_3g_4S.$$
From these relations it follows that the group we computed and the group
computed in \cite{kato} are conjugate.

\end{document}